\documentclass[12pt,reqno]{amsart}

\usepackage{amsmath,amsfonts,amsthm,amssymb,amsxtra, mathrsfs
}
\usepackage{bbm} 
\usepackage{hyperref} 

\newtheorem{theorem}{Theorem}

\newtheorem{corollary}[theorem]{Corollary}

\theoremstyle{definition}

\theoremstyle{remark}

\newcommand\eps\varepsilon

\renewcommand\Re{\mathop{\mathrm{Re}}\nolimits}

\begin{document}
\title[Reverse Heinz type inequality, Mather beta function and coefficient estimates]{Reverse Heinz type inequality, Mather beta function and coefficient estimates}
\author{Vladimir Bo\v zin}
\address{University of Belgrade, Faculty of Mathematics, 11000, Belgrade, Serbia}
\email{vladimir.bozin@matf.bg.ac.rs}
\author{Petar Melentijevi\'c}
\address{University of Belgrade, Faculty of Mathematics, 11000, Belgrade, Serbia}
\email{petar.melentijevic@matf.bg.ac.rs}
\keywords{Heinz estimate, Gaussian curvature, Fourier coefficients, Mather beta function, sharp inequalities, Elliptic integrals}
\subjclass{33E05, 31A05, 53A10}

\begin{abstract}
In this paper, we prove that for two cyclically ordered sequences of complex numbers of unit modulus $(\xi_j)_{j=1}^{n}$ and $(\zeta_j)_{j=1}^{n}$ the inequality:
$$\bigg|\sum_{m=1}^{n}(\xi_{m+1}-\xi_m)\zeta_m\bigg|^2+ \bigg|\sum_{m=1}^{n}(\xi_{m+1}-\xi_m)\zeta_m^{-1}\bigg|^2\leqslant 4n^2\sin^2\frac{\pi}{n}$$
holds for $n\geqslant 4$. As a consequence, for $n=4,$ we disprove Hall's conjecture on Heinz type inequality for harmonic self-mappings of the unit disk. The connection with the theory of elliptic billiards via the Mather beta function will also be given. 
\end{abstract}
\maketitle

\section{Introduction and main results}

\subsection{Motivation for the problem}Let $\Sigma:=\{(u,v, F(u,v)\}$ be a minimal graph in $\mathbb{R}^3 $ over the unit disk ${D_R}:=\{(u,v) \in \mathbb{R}^2: u^2+v^2<R^2\}\subset \mathbb{C}\equiv \mathbb{R}^2.$ Gaussian curvature of this surface is given by
$$K(P)=\frac{F_{uu}F_{vv}-F^2_{uv}}{(1+F_u^2+F_v^2)^2},$$
where $P=(u,v,F(u,v)).$
In 1952, E. Heinz \cite{Heinz} proved the inequality 
\begin{equation}
\label{Gaus}
|K(0)|<\frac{c_0}{R^2}, 
\end{equation}
for some constant $c_0>0.$ This result was further improved by H. Hopf \cite{Hopf}, who proved the estimate
\begin{equation}
\label{Hopf}
(1+|\nabla F(0)|^2)|K(0)|<\frac{c_1}{R^2},  
\end{equation}
with some absolute constant $c_1>0.$ 
Recently, in \cite{KalMel} \eqref{Gaus} is proved with $c_0=\frac{\pi^2}{2},$ while in the paper \cite{KalZhu} stronger inequality \eqref{Hopf} is shown with $c_1=c_0.$ In both papers, the main ingredient in the proofs is a certain comparison theorem of Kalaj, from the unpublished manuscript \cite{Kal}, which generalizes the approach from Finn-Osserman paper \cite{FO} and reduces the problem to estimating the curvature of Scherk-type surfaces.

An interesting approach to finding the best constants in inequalities \eqref{Gaus} and \eqref{Hopf} was given in Duren's monograph \cite{Duren}. Namely, if 
$$f(z)=h(z)+\overline{g(z)}=\sum_{n=0}^{+\infty}{a_n}z^n+\sum_{n=1}^{+\infty}a_{-n}\overline{z}^n=u(z)+iv(z)$$
is a harmonic mapping with dilatation $\frac{g'(z)}{h'(z)}=q(z)^2,$ then the harmonic conformal parameterization of the surface $\Sigma$ is given by
$$u(z)=\Re f(z)=\Re \int_{0}^{z}(h'(\zeta)+g'(\zeta))d\zeta,$$
$$v(z)=\Im f(z)=\Re \int_{0}^{z}i (g'(\zeta)-h'(\zeta))d\zeta,$$
$$T(z)=\Im\int_{0}^{z}2h'(\zeta)q(\zeta)d\zeta.$$
For this surface, the Gaussian curvature is 
$$K=-\frac{4|q'(z)|^2}{|h'(z)|^2(1+|q(z)|^2)^4},$$
while 
$$1+|\nabla F(0)|^2=\bigg(\frac{1+|q(0))^2}{1-|q(0)|^2}\bigg)^2.$$
Using Schwarz lemma and $h'(0)=a_1, g'(0)=a_{-1},$ we get
\begin{align*}
 (1+|\nabla F(0)|^2)|K(0)|&=\frac{4|q'(0)|^2}{|h'(0)|^2(1-|q(0)|^4))^2}\\
 & \leqslant \frac{4(1-|q(0|^2)^2}{|h'(0)|^2(1-|q(0)|^4)^2}\\
 &=\frac{4}{|h'(0)|^2(1+|q(0)|^2)^2}\\
 &=\frac{4}{(|a_1|+|a_{-1}|)^2}.
\end{align*}
A very natural conjecture for these coefficients, as posed in \cite{Duren}, which would lead to the sharp form of the estimate \eqref{Hopf} is
$$|a_1|^2+|a_{-1}|^2 \geqslant \frac{8}{\pi^2}.$$
R. R. Hall \cite{Hall1} proved the appropriate inequality for harmonic self-mappings $f$ of the unit disk $\mathbb{D}$ without the additional assumption that the dilatation is square of an analytic function. More precisely, he proved
$$|a_1|^2+|a_{-1}|^2\geqslant \frac{27}{4\pi^2}$$
and consequently, $|K|\leqslant \frac{16\pi^2}{27}.$
Later on, he proved a slightly refined estimate for Gaussian curvature which showed the crucial role played by the assumption on the dilatation of the harmonic function $f.$\\

\subsection{Poisson extension of a step function, Fourier coefficients}In this paper, we will prove a general inequality for $|a_1|^2+|a_{-1}|^2$ for Poisson extensions of step functions. The good exposition literature for such mappings is \cite{TSS} and \cite{Duren}.\\

Following \cite{Duren}, for a step function $\psi:$
\begin{equation}\psi(e^{\imath t})=\begin{cases}
    e^{i \varphi_1}, \quad t_0=0\leqslant t<t_1,\\
    e^{i \varphi_2}, \quad t_1\leqslant t<t_2,\\
    \dots\\
    e^{i \varphi_n}, \quad t_{n-1}\leqslant t<t_n=2\pi
\end{cases}
\end{equation}
Fourier coefficients of its Poisson integral
$$f(re^{i s})=P[\psi](re^{i s})=\frac{1}{2\pi}\int_{0}^{2\pi}\frac{1-r^2}{1-2r\cos(s-t)+r^2}dt$$
are given by
\begin{align}
 a_l=
 &\frac{1}{2\pi i l}
\sum_{m=1}^{n}
e^{i\varphi_m}(e^{-it_{m-1}}-e^{-i l t_m}) \nonumber \\
=
&\frac{1}{2\pi i l}
\sum_{m=1}^{n}
(e^{i\varphi_{m+1}}-e^{i\varphi_m})e^{-i l t_m} \label{al}
\end{align} 
and 
\begin{align}
a_{-l}=&\frac{1}{2\pi i l}
\sum_{m=1}^{n}e^{i\varphi_m}(e^{i l t_{m}}-e^{i l t_{m-1}})\nonumber \\
=&\frac{1}{2\pi i l}\sum_{m=1}^{n}
(e^{i\varphi_{m}}-e^{i\varphi_{m+1}})e^{i l t_m} \label{a-l},
\end{align}
where $\varphi_{n+1}=\varphi_1.$

We will say the sequence of complex numbers of modulus one is {\it cyclically ordered}, if their arguments can be chosen to be  non decreasing and belonging to an interval of size $2\pi$. 

Our main result reads as follows:
\begin{theorem}
For two sequences of numbers $(\xi_m)_{m=1}^n,(\zeta_m)_{m=1}^n \in \mathbb{T}=\partial\mathbb{D}$, $n \geqslant 4$ and $(\xi_m)$ is cyclically ordered, we have the inequality
$$\bigg|\sum_{m=1}^{n}(\xi_{m+1}-\xi_m)\zeta_m\bigg|^2+ \bigg|\sum_{m=1}^{n}(\xi_{m+1}-\xi_m)\zeta_m^{-1}\bigg|^2\leqslant 4 n^2\sin^2\frac{\pi}{n}.$$
The equality for $n>4$ holds when both sequences are vertices of regular $n$-gons. For $n=4$ equality is achieved whenever there are two members of the first sequence that are opposite for the suitable choice of the second sequence; the LHS value 32 is also the maximal value of the LHS sum for $n=3$, with equality achieved under the same condition that two members of the first sequence are opposite.
\end{theorem}
Note that, by the monotonicity of our bound with $n$, and the fact that if our cyclical sequence has two repeated points the LHS expression reduces to the expression for $n-1$, we may assume that the points of our cyclically ordered sequence are different.

An immediate consequence is the following:
\begin{corollary}
Let $f=h+\overline{g}$ be a harmonic self-mapping of the unit disk with the dilatation $\frac{g'}{h'}$ equal to the square of the Blaschke product $B(z)$  and $a_1, a_{-1}$ are its Fourier coefficients,
then
$$|a_1|^2+|a_{-1}|^2\leqslant \frac{n^2}{\pi^2} \sin^2\frac{\pi}{n}$$
Specifically, for $n=4,$ we get
$$|a_1|^2+|a_{-1}|^2\leqslant \frac{8}{\pi^2}$$
for harmonic mapping whose dilatation is the square of a M\"obius transform.
\end{corollary}

Note that the partial summation formulas and reversing order by conjugation will give us the same inequality if either of the two sequences is cyclically ordered. 

Our proof of Theorem 1 uses results from the theory of elliptic billiards. Namely, we use Poncelet's Porism for polygons of maximal perimeter, i.e., $n$-gons of maximal perimeter inscribed in an ellipse form a continuous $1$-parameter family sharing a common inscribed confocal ellipse. These $n$-gons can be understood as a closed path of an elliptic billiard with rotation number $\rho=1/n$.

\begin{theorem}[see \cite{Berger87}, Theorem 17.6.6]
\label{thm:berger_max_perimeter}
Let $C$ be an ellipse and $n \geqslant 3$ an integer. Among the convex polygons with $n$ distinct vertices inscribed in $C$, there exist some with maximum perimeter. In fact, there are infinitely many such maximum-perimeter polygons (MPP); one vertex of an MPP can be chosen arbitrarily on $C$. Furthermore, the sides of all $n$-vertex MPP's are tangent to the same ellipse $C'_n$, homofocal with $C$.
\end{theorem}

Let $E$ be an ellipse with semi-axes $a > b > 0$, with $c = \sqrt{a^2 - b^2}$. The family of confocal caustics $E_\lambda$ is parameterized by $\lambda \in (0, b^2)$ as
\[
\frac{x^2}{a^2 - \lambda} + \frac{y^2}{b^2 - \lambda} = 1.
\]
Each caustic $E_\lambda$ has eccentricity $f = \frac{c}{\sqrt{a^2 - \lambda}}$ and perimeter $|E_\lambda|$, and there is a unique correspondence of $\lambda$ to a rotation number $\rho$. A recent paper by Bialy \cite{Bially} expresses the Mather $\beta$ function, measuring the average chord length of one segment of a billiard trajectory of rotation number $\rho$ explicitly in terms of elliptic functions. Note that the maximal perimeter is then simply $n \beta(1/n)$.

\begin{theorem}[see \cite{Bially}, Corollary 2.2]
\label{bialy_mather}
The Mather $\beta$-function $\beta(\rho)$ for the ellipse $E$ satisfies:
\[
\beta(\rho) = \frac{2a\sqrt{\lambda}}{b} - 2\sqrt{a^2 - \lambda} \, E(\phi, k) + \rho |E_\lambda|,
\]
where
\[
\phi = \arcsin \frac{\sqrt{\lambda}}{b}, \quad k = 1/f,
\]
and $E(\phi, k)$ is the incomplete elliptic integral of the second kind.
\end{theorem}

\section{Proof of the main result}

From \eqref{al} and \eqref{a-l}, we have:
\[
|a_1|^2+|a_{-1}|^2
=
\frac{1}{4\pi^2}
\left(
\left|
\sum_{m=1}^{n}
\left(e^{i\varphi_{m+1}}-e^{i\varphi_m}\right)e^{it_m}
\right|^2
+
\left|
\sum_{m=1}^{n}
\left(e^{i\varphi_{m+1}}-e^{i\varphi_m}\right)e^{-it_m}
\right|^2
\right).
\]
Using the polarization identity $|a+b|^2+|a-b|^2=2(|a|^2+|b|^2)$ for
$a=\sum_{m=1}^{n}(e^{i\varphi_{m+1}}-e^{i\varphi_m})\cos t_m$ and $b=i\sum_{m=1}^{n}(e^{i\varphi_{m+1}}-e^{i\varphi_m})\sin t_m,$
this is further equal to:
\[
\frac{1}{2\pi^2}
\left(
\left|
\sum_{m=1}^{n}
\left(e^{i\varphi_{m+1}}-e^{i\varphi_m}\right)\cos t_m
\right|^2
+
\left|
\sum_{m=1}^{n}
\left(e^{i\varphi_{m+1}}-e^{i\varphi_m}\right)\sin t_m
\right|^2
\right),
\]
Note that the expression in brackets can be interpreted as the norm of the quaternion
\[
A=
\sum_{m=1}^{n}
\left(e^{i\varphi_{m+1}}-e^{i\varphi_m}\right)e^{jt_m},
\]
where $i, j, k=ij$ are the imaginary units. 
We have
\[
A
=
2i
\sum_{m=1}^{n}
\sin\frac{\varphi_{m+1}-\varphi_m}{2}
\left(
\cos\frac{\varphi_{m+1}+\varphi_m}{2}
+
i\sin\frac{\varphi_{m+1}+\varphi_m}{2}
\right)
\left(
\cos t_m+j\sin t_m
\right)
\]
and its expanded form is
\[
\begin{aligned}
A
&=
-2i
\sum_{m=1}^{n}
\sin\frac{\varphi_{m+1}-\varphi_m}{2}
\Biggl[
\cos\frac{\varphi_{m+1}+\varphi_m}{2}\cos t_m
\\
&\qquad
+i\sin\frac{\varphi_{m+1}+\varphi_m}{2}\cos t_m
\\
&\qquad
+j\cos\frac{\varphi_{m+1}+\varphi_m}{2}\sin t_m
\\
&\qquad
+k\sin\frac{\varphi_{m+1}+\varphi_m}{2}\sin t_m
\Biggr]
\\
&=
i\left(
S_1+iS_2+jS_3+kS_4
\right).
\end{aligned}
\]

Since $\|A\|=\max_{\|Q\|=1}|\langle A, Q\rangle|$ and every quaternion of the unit norm can be represented as a product $e^{ui}e^{kt}e^{jv},$ for some $t, u, v\in \mathbb{R},$ we get
\[
\|A\|
=
\|S_1+iS_2+jS_3+kS_4\|
=
\sup_{\substack{t\in[0,2\pi)\\u,v\in\mathbb{R}}}
\left|
\left\langle
A,
e^{iu}e^{kt}e^{jv}
\right\rangle
\right|,
\]
and hence:
\begin{align*}
\sup_{\substack{
0<\varphi_m<\varphi_{m+1}<2\pi\\
0\le t_1,t_2,\ldots,t_n\le2\pi}}
\|A\|
=&
\sup_{\varphi_m,t_m,t,u,v}
\left|
\left\langle
\sum_{m=1}^{n}
\left(e^{i\varphi_{m+1}}-e^{i\varphi_m}\right)e^{jt_m},
e^{iu}e^{kt}e^{jv}
\right\rangle
\right|\\
=&
\sup_{\varphi_m,t_m,t,u,v}
\left|
\left\langle
e^{-iu}
\sum_{m=1}^{n}
\left(e^{i\varphi_{m+1}}-e^{i\varphi_m}\right)
e^{jt_m}e^{jv},
e^{kt}
\right\rangle
\right|\\
=&
\sup_{\varphi_m,t_m,t,u,v}
\left|
\left\langle
\sum_{m=1}^{n}
\left(
e^{i(\varphi_{m+1}-u)}
-
e^{i(\varphi_m-u)}
\right)
e^{j(t_m+v)},
e^{kt}
\right\rangle
\right|\\
=&
\sup_{\substack{t\in[0,2\pi)\\ \varphi_m,t_m}}
\left|
\left\langle
\sum_{m=1}^{n}
\left(e^{i\varphi_{m+1}}-e^{i\varphi_m}\right)e^{jt_m},
e^{kt}
\right\rangle
\right|\\
=&
\sup_{\substack{t\in[0,2\pi)\\ \varphi_m,t_m}}
\left|
\left\langle
S_1+iS_2+jS_3+kS_4,
\cos t+k\sin t
\right\rangle
\right|
\end{align*}
Therefore,
$$\sup_{\substack{
0<\varphi_m<\varphi_{m+1}<2\pi\\
0\le t_1,t_2,\ldots,t_n\le2\pi}}
\|A\|
=
\sup_t
\left|
S_1\cos t+S_4\sin t
\right|
=
\sup_t
\left(
S_1\cos t+S_4\sin t
\right).
$$

Now, our problem reduces to finding the maximum of the function:
\begin{align}
&\Phi(\varphi_1,\dots,\varphi_n,t_1,\dots,t_n,t)\nonumber \\
=&
2\sum_{m=1}^{n}
\sin\frac{\varphi_{m+1}-\varphi_m}{2}
\left(
\cos\frac{\varphi_{m+1}+\varphi_m}{2}
\cos t_m\cos t
+
\sin\frac{\varphi_{m+1}+\varphi_m}{2}
\sin t_m\sin t
\right). \label{Phi}
\end{align}

We will first determine the values of \(t_m\) so that the last function reaches its maximal value. From the system 
$$
\frac{\partial\Phi}{\partial t_m}
=
2\sin\frac{\varphi_{m+1}-\varphi_m}{2}
\left(
-
\cos\frac{\varphi_{m+1}+\varphi_m}{2}
\sin t_m\cos t
+
\sin\frac{\varphi_{m+1}+\varphi_m}{2}
\cos t_m\sin t
\right)
=0,
$$
we infer, using the assumption that the points of our cyclically ordered sequence are different and $\sin {2t}\neq 0$, for stationary points we have

$$
\tan t_m
=
\frac{\sin t_m}{\cos t_m}
=
\frac{
\sin\frac{\varphi_{m+1}+\varphi_m}{2}\sin t
}{
\cos\frac{\varphi_{m+1}+\varphi_m}{2}\cos t
}
=
\tan t\,
\tan\frac{\varphi_{m+1}+\varphi_m}{2}.
$$
This implies
\[
\sin t_m
=
\pm
\frac{
\sin t\,
\sin\frac{\varphi_{m+1}+\varphi_m}{2}
}{
\sqrt{
\sin^2 t\,
\sin^2\frac{\varphi_{m+1}+\varphi_m}{2}
+
\cos^2 t\,
\cos^2\frac{\varphi_{m+1}+\varphi_m}{2}
}
}
\]

\[
\cos t_m
=
\pm
\frac{
\cos t\,
\cos\frac{\varphi_{m+1}+\varphi_m}{2}
}{
\sqrt{
\sin^2 t\,
\sin^2\frac{\varphi_{m+1}+\varphi_m}{2}
+
\cos^2 t\,
\cos^2\frac{\varphi_{m+1}+\varphi_m}{2}
}
}.
\]

Note that case $\sin 2t=0$ can be examined separately, and it is easily seen that then the maximal value of expression under this condition is reduced to projections along either the horizontal or vertical axis, with the maximum of $\Phi$ equal to $4$ achieved when there are two members of the first sequence that are opposite.
Namely, for $t=0,$ we get:
\begin{align*}
\Phi(\varphi_1,\dots, \varphi_n, t_1, t_2,\dots,t_n,0)&
=\sum_{m=1}^{n}(\cos\varphi_m-\cos\varphi_{m+1})\cos t_m\\
&\leqslant \sum_{m=1}^{n}|\cos\varphi_m-\cos\varphi_{m+1}|\leqslant 4,
\end{align*}
since the variation of the cosine function in $[0,2\pi]$ is equal to 4. Similarly we get the analogous estimate for $t=\pi.$
This value will correspond to the maximum in the case $n$ equal to $3$ or $4,$ since the value of $\Phi$ in stationary points for $\sin 2t\neq 0$ is $\sqrt{2}n\sin\frac{\pi}{n}$ (we will prove this in the following lines) which is smaller than $4$ for $n=3$ and equal to $4$ for $n=4.$

Otherwise, the value of the function $\Phi$ for critical values of $t_m$ is
\[
\begin{aligned}
&\Phi(\varphi_1,\ldots,\varphi_n,t_1,\ldots,t_n,t)\\
&=
2
\sum_{m=1}^{n}
\sin\frac{\varphi_{m+1}-\varphi_m}{2}
\left(
\pm
\frac{
\cos^2 t\,
\cos^2\frac{\varphi_{m+1}+\varphi_m}{2}
+
\sin^2 t\
\sin^2\frac{\varphi_{m+1}+\varphi_m}{2}
}{
\sqrt{
\sin^2 t\,
\sin^2\frac{\varphi_{m+1}+\varphi_m}{2}
+
\cos^2 t\,
\cos^2\frac{\varphi_{m+1}+\varphi_m}{2}
}
}
\right)
\\
&=
\pm
2
\sum_{m=1}^{n}
\sin\frac{\varphi_{m+1}-\varphi_m}{2}
\sqrt{
\sin^2 t\,
\sin^2\frac{\varphi_{m+1}+\varphi_m}{2}
+
\cos^2 t\,
\cos^2\frac{\varphi_{m+1}+\varphi_m}{2}
}\\
&=\pm \sum_{m=1}^{n}\sqrt{\sin^2t(\cos\varphi_m-\cos\varphi_{m+1})^2+\cos^2t(\sin\varphi_m-\sin\varphi_{m+1})^2}
\end{aligned}
\]
and, clearly, the plus signs correspond to the maximum value, achievable for the suitable choice of $t_m$, since the choice of these signs is at our disposal. Denote $x=\sin^2 t$ and consider 
\[
f_n(x)=\sup_{\varphi_1<\varphi_2<\dots<\varphi_n}
\sum_{m=1}^{n}
\sqrt{x(\cos\varphi_m-\cos\varphi_{m+1})^2+(1-x)(\sin\varphi_m-\sin\varphi_{m+1})^2
}.
\]
where $\varphi_n<\varphi_{n+1}=\varphi_1+2\pi.$
\\
Now note that since we assume that the sequence $e^{i \varphi_m}$ is a cyclically ordered sequence of distinct points, we can interpret this as the maximal perimeter of the $n$-gon, inscribed in an ellipse with semi-axes $\sqrt x$ and $\sqrt{1-x}$. 

By Poncelet's Porism, the maximal perimeter for fixed $x$ is given by the corresponding value of the Mather $\beta$ function for $\rho=1/n$. Fix this $\rho$, and consider this as a function of $x$. We will show that for $n > 4$ this function achieves a maximum at $x=1/2$, and for $n=3$, it will have a minimum there, while for $n=4$ the function will actually be a constant. We will compute the derivative with respect to $x$ for $x<1/2$ and show that the function is increasing for $n>4$.

We proceed with our computation, using well-known properties of elliptic functions (see, e.g., \cite{Akhiezer90}). 
Let us start from the following expression for the Mather beta function:
$$\beta_{\rho}=2\sqrt{\frac{\lambda(1-x)}{x}}-2\sqrt{1-x-\lambda}E(\varphi,k)+4\sqrt{1-x-\lambda}\rho E_c(k),$$
where $\rho=\frac{F(k)}{2K(k)}.$
Let  $D = \frac{d}{dx} \rho\vert_{g=\text{const.}}$ and denote
$S(x,\lambda) = \sqrt{1-x-\lambda}$, $m = \kappa^2$, $\Delta = \sqrt{1-\kappa^2 \sin^2\varphi} = \sqrt{1-m\sin^2\varphi}$, 
$E = E(\varphi, m)$, $E_c = E_c(m)$, $\frac{F}{K} = 2g = \text{const.}$, $Z = E - \frac{F}{K} E_c$.

Since $\frac{\partial E}{\partial m} = \frac{E-F}{2m}$, $\frac{\partial E_c}{\partial m} = \frac{E_c-K}{2m}$ and 
$\beta = 2\sqrt{\frac{\lambda(1-x)}{x}} - 2SZ$, we get 
$$ D\beta = D\left(2\sqrt{\frac{\lambda(1-x)}{x}}\right) - 2(DS)Z - 2SDZ,$$
while
\begin{align*}
DZ &= DE - \frac{F}{K} DE_c \\
&= \frac{\partial E}{\partial \varphi} D\varphi + \frac{\partial E}{\partial m} Dm - \frac{F}{K} \frac{\partial E_c}{\partial m} Dm \\
&= \Delta D\varphi + \frac{E-F-\frac{F}{K}(E_c-K)}{2m} Dm \\
&= \Delta D\varphi + \frac{Z}{2} \frac{Dm}{m}.
\end{align*}
Further:
\begin{align*}
D\beta &= D\left(2\sqrt{\frac{\lambda(1-x)}{x}}\right) - 2(DS)Z - 2S \cdot DZ \\
&= D\left(2\sqrt{\frac{\lambda(1-x)}{x}}\right) - 2(DS)Z - 2S \left(\Delta D\varphi + Z \frac{Dm}{2m}\right).
\end{align*}

From $m = \frac{1-2x}{S^2}$ we see that 
$$\frac{Dm}{m} = \frac{-2}{1-2x} - \frac{2DS}{S},$$
so, we have:
\begin{align*}
D\beta &= D\left(2\sqrt{\frac{\lambda(1-x)}{x}}\right) - 2S\Delta D\varphi - 2DS \cdot Z - S \cdot Z \cdot \left(\frac{-2}{1-2x} - \frac{2DS}{S}\right) \\
&= D\left(2\sqrt{\frac{\lambda(1-x)}{x}}\right) - 2S\Delta D\varphi + \frac{2SZ}{1-2x} \\
&= \frac{\partial}{\partial x}\left(2\sqrt{\frac{\lambda(1-x)}{x}}\right) - 2S\Delta \frac{\partial\varphi}{\partial x} \\
&\quad + \left[ \frac{\partial}{\partial \lambda}\left(2\sqrt{\frac{\lambda(1-x)}{x}}\right) - 2S\Delta \frac{\partial\varphi}{\partial \lambda} \right] D\lambda + \frac{2\sqrt{1-x-\lambda}}{1-2x} Z.
\end{align*}

Since $\frac{\partial}{\partial x}\left(2\sqrt{\frac{\lambda(1-x)}{x}}\right) = -\frac{\sqrt{\lambda}}{x\sqrt{x(1-x)}}$, 
$\frac{\partial}{\partial \lambda}\left(2\sqrt{\frac{\lambda(1-x)}{x}}\right) = \frac{1}{\sqrt{\lambda}}\sqrt{\frac{1-x}{x}}$
and
\begin{align*}
S \cdot \Delta &= \sqrt{1-x-\lambda} \cdot \sqrt{1-\kappa^2\sin^2\varphi} \\
&= \sqrt{1-x-\lambda} \cdot \sqrt{1 - \frac{1-2x}{1-x-\lambda} \cdot \frac{\lambda}{x}} \\
&= \sqrt{1-x-\lambda} \cdot \frac{\sqrt{x(1-x-\lambda) - \lambda(1-2x)}}{\sqrt{x(1-x-\lambda)}} \\
&= \frac{\sqrt{1-x} \cdot \sqrt{x-\lambda}}{\sqrt{x}},
\end{align*}
while  $\frac{\partial\varphi}{\partial\lambda}= \frac{1}{2\sqrt{\lambda}\sqrt{x-\lambda}}$ and $\frac{\partial\varphi}{\partial x}= -\frac{\sqrt{\lambda}}{2x\sqrt{x-\lambda}},$
we easily see that the expression with $D\lambda$ is equal to:
\begin{align*}
\frac{1}{\sqrt{\lambda}}\sqrt{\frac{1-x}{x}} - 2 \frac{\sqrt{1-x}\sqrt{x-\lambda}}{\sqrt{x}} \cdot \frac{1}{2\sqrt{\lambda}\sqrt{x-\lambda}} &= \sqrt{\frac{1-x}{x}}\frac{1}{\sqrt{\lambda}} - \sqrt{\frac{1-x}{x}}\frac{1}{\sqrt{\lambda}} = 0!
\end{align*}

Therefore:
\begin{align*}
D\beta &= -\frac{\sqrt{\lambda}}{x\sqrt{x(1-x)}} - 2 \frac{\sqrt{1-x}\sqrt{x-\lambda}}{\sqrt{x}} \cdot \left(-\frac{\sqrt{\lambda}}{2x\sqrt{x-\lambda}}\right) + 2Z\frac{\sqrt{1-x-\lambda}}{1-2x} \\
&= -\frac{\sqrt{\lambda}}{x\sqrt{x}} \left( \sqrt{1-x} - \frac{1}{\sqrt{1-x}} \right) + \frac{2\sqrt{1-x-\lambda}}{1-2x} Z \\
&= -\frac{\sqrt{\lambda}}{\sqrt{x(1-x)}} + \frac{2\sqrt{1-x-\lambda}}{1-2x} Z, 
\end{align*}
and $\beta$ is increasing as a function on $\varphi$ if and only if
$$Z \geqslant \frac{\kappa^2 \sin\varphi \cos\varphi}{2\sqrt{1-\kappa^2\sin^2\varphi}}! $$

Consider the function 
$$L(\varphi)=\int_{0}^{\varphi}\sqrt{1-k^2\sin^2\theta}d\theta-\frac{E_c(k)}{K(k)}\int_{0}^{\varphi}\frac{1}{\sqrt{1-k^2\sin^2\theta}}d\theta-\frac{k^2\sin\varphi\cos\varphi}{2\sqrt{1-k^2\sin^2\varphi}}.$$
Since
$$L'(\varphi)=\sqrt{1-k^2\sin^2\varphi}-\frac{E_c(k)}{K(k)}\frac{1}{\sqrt{1-k^2\sin^2\varphi}}-\frac{k^2(k^2\sin^4\varphi+1-2\sin^2\varphi)}{2(1-k^2\sin^2\varphi)^{\frac{3}{2}}},$$
we get
$$L'(\varphi)(1-k^2\sin^2\varphi)=\frac{1}{2}(1-k^2\sin^2\varphi)^2-\frac{E_c(k)}{K(k)}(1-k^2\sin^2\varphi)+\frac{1}{2}(1-k^2)=:R(\varphi).$$
Note that $R(0)=1-\frac{1}{2}k^2-\frac{E_c(k)}{K(k)}>0$ and $R(\frac{\pi}{2})=(1-k^2)(1-\frac{1}{2}k^2-\frac{E_c(k)}{K(k)})>0,$
while $R\big(\arcsin\frac{1-\frac{E_c(k)}{K(k)}}{k^2}\big)=\frac{1}{2}(1-k^2-\frac{E_c^2(k)}{K^2(k)})<0$ since $\sqrt{1-k^2}<\frac{E_c(k)}{K(k)}<1-\frac{1}{2}k^2$ for $k \in (0,1).$ This means that there exist some $\varphi_0(k)$ and $\varphi_1(k)$ such that $L$ increases for $\varphi$ from $0$ to $\varphi_0(k)$ then decreases from $\varphi_0(k)$ to $\varphi_1(k)$ and increases for $\varphi$ from $\varphi_1(k)$ to $\frac{\pi}{2}.$ Therefore, $L$ is positive on $(0,\overline{\varphi_k})$ and negative on $(\overline{\varphi}_k,1).$
Let us prove that $F(\overline{\varphi_k},k)=\frac{1}{2}K(k).$ \\

Introducing the change of variable $\sin t=\frac{\cos\theta }{\sqrt{1-k^2\sin^2\theta}}$  in the integral for $F$ we get $K(k)=F(\varphi_1,k)+F(\varphi_2,k)$ for $\sin\varphi_2=\frac{\cos\varphi_1}{\sqrt{1-k^2\sin^2\varphi_1}}, $i. e. $1-k^2\sin^2\varphi_1\sin^2\varphi_2=\cos^2\varphi_1+\cos^2\varphi_2.$ For $\varphi_1=\varphi_2=\overline{\varphi}_k$ with $\sin\overline{\varphi}_k=\frac{\cos\overline{\varphi}_k}{\sqrt{1-k^2\sin^2\overline{\varphi}_k}},$ we have $2F(\overline{\varphi_k},k)=K(k),$ which is to be shown. 

The function $Z$ also satisfies a similar addition formula $Z(\psi,k)=Z(\varphi_1,k)+Z(\varphi_2,k)-k^2\sin\varphi_1\sin\varphi_2\sin\psi.$ Since $0=Z(\frac{\pi}{2},k)=2Z(\overline{\varphi}_k,k)-k^2\sin^2\overline{\varphi}_k,$ we finally get $Z(\overline{\varphi}_k,k)=\frac{1}{2}k^2\sin^2\overline{\varphi}_k=\frac{k^2\sin\overline{\varphi}_k\cos\overline{\varphi}_k}{2\sqrt{1-k^2\sin^2\overline{\varphi}_k}}$ and $L(\varphi)\geqslant 0$ for $\varphi \in (0,1).$
Hence, for $n>4,$ $f_n(x)$ has its maximum in $\frac{1}{2}$ equal to
$\sqrt{2}\sup\sum_{m=1}^n\sin\frac{\varphi_{m+1}-\varphi_m}{2},$ which by concavity of the sine function is equal to $\sqrt{2}n\sin\frac{\pi}{n}.$

\section{Acknowledgements} The first author is partially supported by MPNTR grant no. 174032, and the second author is partially supported by MPNTR grant no. 174017.\\


\bibliographystyle{amsalpha}

\end{document}